\documentclass[12pt]{article}

\usepackage{amsmath, amssymb}

\newtheorem{thm}{Theorem}
\newtheorem{prop}{Proposition}
\newtheorem{cor}{Corollary}
\newtheorem{lem}{Lemma}

\title{Ample line bundles \\ on  certain toric fibered 3-folds}
\author{Shoetsu Ogata}

\begin{document}
\maketitle

\begin{abstract}
Let $X$ be a projective nonsingular toric 3-fold 
with a surjective torus equivariant morphism onto the projective line 
or a nonsingular toric surface not isomorphic to the projective plane.
Then we prove that an ample line bundle on $X$ is always normally
generated.
\end{abstract}
\section*{Introduction}

Let $X$ be a projective algebraic variety and let $\cal{L}$ an ample line
bundle on it.  If the multiplication map 
\begin{equation}\label{0;0}
\mbox{Sym}^k \Gamma(X, {\cal L}) \longrightarrow \Gamma(X, {\cal{L}}^{\otimes k})
\end{equation}
of the $k$-th symmetric power of the
global sections of $\cal{L}$ to the global sections of the $k$-th tensor
product is surjective for all positive integers $k$, then
Mumford \cite{Mf} calls $\cal{L}$ {\it normally generated}.
A normally generated ample line bundle is always  very ample, but not conversely.

If $X$ is a toric variety of dimension $n$ and $\cal{L}$ is an ample line bundle on it, then 
Ewald and Wessels \cite{EW} showed that ${\cal L}^{\otimes k}$ is very ample for $k\ge n-1$,
and Nakagawa \cite{N} showed that the multiplication map
\begin{equation*}
\Gamma(X, {\cal L}^{\otimes k})\otimes \Gamma(X,{\cal L})
\longrightarrow \Gamma(X, {\cal L}^{\otimes (k+1)})
\end{equation*}
is surjective for $k\ge n-1$.
We  know that there exists a polarized toric variety $(X, {\cal L})$ of dimension
$n\ge3$ such that ${\cal L}^{\otimes (n-2)}$ is not normally generated.
We also know that any ample line bundle on a nonsingular toric variety is
always very ample (see \cite[Corollary 2.15]{Od}).
Ogata \cite{Og} showed that an ample line bundle $L$ on a nonsingular toric
3-fold $X$  is normally generated if 
 the adjoint bundle ${\cal L}+K_X$ is not big.

If a toric variety $X$ of dimension $n\ge2$ has a surjective torus equivariant
morphism $\varphi: X \to Y$ onto a toric variety $Y$ of dimension $r$ ($1\le r<n$)
with connected fibers, then we call $X$ a {\it toric fibered $n$-fold over $Y$}.

In this paper we restrict $X$ to a certain class of toric fibered 3-folds.
\begin{thm}\label{0:t1}
Let $X$ be a nonsingular projective toric fibered 3-fold over the projective line.
Then an ample line bundle on $X$ is always normally generated.
\end{thm}

Since a nonsingular toric surface not isomorphic to $\mathbb{P}^2$ has
a  torus equivariant surjective morphism onto $\mathbb{P}^1$,
Theorem~\ref{0:t1} implies the following corollary.

\begin{cor}
Let $X$ be a nonsingular projective toric fibered 3-fold over a nonsingular
toric surface not isomorphic to the projective plane.
Then an ample line bundle on $X$ is always normally generated.
\end{cor}

On a toric variety $X$ of dimension $n$,
the space 
$\Gamma(X, \cal{L})$ of global sections of an 
ample line bundle $\cal{L}$ is parametrized by lattice points in a lattice polytope $P$
of dimension $n$ (see, for instance, Oda's book\cite[Section 2.2]{Od} 
or Fulton's book\cite[Section 3.5]{Fu}). 
If $X$ has a surjective morphism $\varphi: X\to \mathbb{P}^1$ onto the projective
line, then the corresponding polytope $P$ has a special shape.
From this fact, we shall prove Theorem~\ref{0:t1}.

In Section 3, we prove the same statement of Theorem~\ref{0:t1} under the
assumption that one invariant fiber of $\varphi$ is irreducible.  This is given
as Proposition~\ref{3;p1}.
Full statement is proved in Section 4 as Proposition~\ref{4;p1}.
In the end of this paper, we remark that nonsingularity condition is necessary
by giving an example.

\bigskip

\section{Toric varieties and lattice polytopes}

In this section we recall the fact about toric varieties and ample line bundles on them
and corresponding lattice polytopes from Oda's book \cite{Od} or Fulton's book \cite{Fu}.

Let $N\cong \mathbb{Z}^n$ be a free abelian group of rank $n$ and $M:=
\mbox{Hom}(N, \mathbb{Z})$ its dual with the pairing $\langle\cdot, \cdot\rangle:
M\times N \to \mathbb{Z}$.
By scalar extension to real numbers $\mathbb{R}$, we have real vector spaces
$N_{\mathbb{R}}:= N\otimes_{\mathbb{Z}}\mathbb{R}$ and $M_{\mathbb{R}}:=
M\otimes_{\mathbb{Z}}\mathbb{R}$.
We also have the paring of $M_{\mathbb{R}}$ and $N_{\mathbb{R}}$ by  scalar extension,
which is denoted by the same symbol $\langle\cdot, \cdot\rangle$. 

The group ring $\mathbb{C}[M]$ defines an algebraic torus $T_N:= \mbox{Spec}
\mathbb{C}[M]\cong (\mathbb{C}^{\times})^n$ of dimension $n$.
Then the character group $\mbox{Hom}_{\rm gr}(T_N, \mathbb{C}^{\times})$ of
the algebraic torus $T_N$ coincides with $M$.
For $m\in M$ we denote the corresponding character by ${\bf e}(m): T_N \to
\mathbb{C}^{\times}$.

Let $\Delta$ be a finite complete fan of $N$ and $X(\Delta)$  denote the toric variety
defined by $\Delta$.
Set $N_0:=\mathbb{Z}$ and $\Delta_0:=\{\mathbb{R}_{\le0}, \{0\}, \mathbb{R}_{\ge0}\}$.
Then we have $X(\Delta_0)=\mathbb{P}^1$.
If a surjective morphism $\varphi: X(\Delta) \to \mathbb{P}^1$ is torus equivariant, 
then it defines a morphism of fans 
$\varphi^{\sharp}: (N, \Delta) \to (N_0, \Delta_0)$.
Moreover, if fibers of 
$\varphi$ are connected, then $\varphi^{\sharp}(N)=N_0$.
Set $N_0^{\vee}$ the dual to $N_0$.  Then the dual homomorphism 
$\varphi^*: N_0^{\vee} \to M
=N^{\vee}$ maps $N_0^{\vee}$ as a saturated submodule in $M$.
Thus we have a direct sum decomposition $M\cong M'\oplus N_0^{\vee}$.

Set $N_f:=(\varphi^{\sharp})^{-1}(0)$ and 
$\Delta_f:=\{\sigma\in \Delta; \varphi^{\sharp}(\sigma)=0\}$.
Then 
$\Delta_f$ is a fan of $N_f$ and the toric variety $X(\Delta_f)$ is 
a general fiber of 
$\varphi: X(\Delta) \to \mathbb{P}^1$.

We define a {\it lattice polytope} as the convex hull $P:= \mbox{Conv}\{m_1,
\dots, m_r\}$ of a finite subset $\{m_1, \dots, m_r\}$ of $M$ in $M_{\mathbb{R}}$.
We define the dimension of a lattice polytope $P$ as that of the smallest
affine subspace $\mathbb{R}(P)$ containing $P$.

Let $X$ be a projective toric variety of dimension $n$ and $\cal{L}$ an ample line bundle
on $X$.   Then there exists a lattice polytope $P$ of dimension $n$ such that
the space of global sections of $\cal{L}$ is described by
\begin{equation}\label{0:1}
\Gamma(X, {\cal L})\cong \bigoplus_{m\in P\cap M}\mathbb{C}{\bf e}(m),
\end{equation}
where ${\bf e}(m)$ is considered as a rational function on $X$ since
$T_N$ is identified with the dense open subset (see
 \cite[Section 2.2]{Od} or \cite[Section 3.5]{Fu}).
Conversely, a lattice polytope $P$ in $M_{\mathbb{R}}$ of dimension $n$
defines a polarized toric variety $(X, L)$ satisfying the equality (\ref{0:1})
(see \cite[Chapter 2]{Od} or \cite[Section 1.5]{Fu}).

The $k$-th  tensor product ${\cal L}^{\otimes k}$ of $\cal{L}$ 
 corresponds to the $k$-th multiple $kP$ of $P$ for a positive integer $k$.
The condition that the multiplication map 
\begin{equation*}
\Gamma(X, {\cal L}^{\otimes k})\otimes \Gamma(X,{\cal L})
\longrightarrow \Gamma(X, {\cal L}^{\otimes (k+1)})
\end{equation*}
is surjective is equivalent to the equality
\begin{equation*}\label{0:3}
(k P)\cap M +(P\cap M) =((k+1)P)\cap M.
\end{equation*}

%\begin{dfn}
A lattice polytope $P$ in $M_{\mathbb{R}}$ is called {\it normal}
if the equality
\begin{equation}\label{0:2}
\overbrace{(P\cap M)+\dots + (P\cap M)}^{\mbox{$k$ times}} = (kP)\cap M
\end{equation}
holds for all positive integers $k$.
%\end{dfn}
This is equivalent to the condition that the equality
\begin{equation}\label{0:3}
(kP)\cap M + P\cap M = ((k+1)P)\cap M
\end{equation}
holds for all positive integers $k$.
We note that an ample line bundle $L$ on a toric variety is normally generated 
if and only if the corresponding lattice polytope $P$ is normal.
We also note that the equality (\ref{0:2}) holds if and only if
for each lattice point 
$v\in (kP)\cap M$, there exists just $k$ lattice points $u_1,  \dots, u_k$ in $P\cap M$ with 
$v=u_1+\dots +u_k$.

In order to prove the normality of  lattice polytopes, the following theorem is useful.
\begin{thm}[Nakagawa \cite{N}]\label{1:t1}
Let $P$ be a lattice polytope in $M_{\mathbb{R}}$ of dimension $n$.
Then we have the equality 
$$
(kP)\cap M +P\cap M =((k+1)P)\cap M
$$
for integer $k\ge n-1$.
\end{thm}

From Theorem~\ref{1:t1} we see that 
 for the normality of $P$ of dimension three, it is enough to show the equality
$$
(P\cap M) +(P\cap M) = (2P)\cap M.
$$

For a face of a lattice polytope, it is called an {\it edge} if it is of dimension one and 
a {\it facet} if of codimension one.
A lattice polytope $P$ of dimension $n$ is called {\it simple} if at each vertex $v$, 
just $n$ edges meet, that is, the convex cone 
$C_v(P):=\mathbb{R}_{\ge0}(P-v)$ is written as 
$$
\mathbb{R}_{\ge0}m_1 + \dots + \mathbb{R}_{\ge0}m_n
$$
with $m_1, \dots, m_n\in M$.
Moreover, if the set 
$\{m_1, \dots, m_n\}$ is a $\mathbb{Z}$-basis of $M$, then 
$P$ is called {\it nonsingular}.  For a face $F$ of $P$, we call
$F$ is {\it nonsingular} if it is nonsingular with respect to the sublattice 
$\mathbb{R}(F)\cap M$.
We note that a face of a nonsingular lattice polytope is also nonsingular.

\bigskip

\section{Polygonal prisms}

For two lattice polytopes $P$ and $Q$ in $M_{\mathbb{R}}$, we define the Minkowski sum as
$$
P+Q:=\{x+y\in M_{\mathbb{R}}; x\in P\ \mbox{and}\ y\in Q\}.
$$
Then $P+Q$ is also a lattice polytope.

\begin{figure}[h]
 \begin{center}
 \begin{tabular}{cc}
 \setlength{\unitlength}{1mm}
  \begin{picture}(50,60)(5,0)
  \put(10,10){\line(1,0){10}}
  \put(10,10){\line(2,1){15}}
  \put(10,10){\line(1,3){10}}
  \put(20,10){\line(2,1){15}}
  \put(25,17.5){\line(1,0){10}}
  \put(35,47.5){\line(1,0){10}}
  \put(20,40){\line(1,0){10}}
  \put(20,40){\line(2,1){15}}
  \put(30,40){\line(2,1){15}}
  \put(20,10){\line(1,3){10}}
  \put(25,17.5){\line(1,3){10}}
  \put(35,17.5){\line(1,3){10}}
  \put(3,2){\makebox(10,10){$0$}}
  \put(10,0){\makebox(10,10){$J_1$}}
  \put(12,14){\makebox(15,10){$J_2$}}
  \put(5,20){\makebox(15,10){$I$}}
  \put(10,10){\circle*{1}}
  \put(20,10){\circle*{1}}
  \put(20,40){\circle*{1}}
  \put(25,17.5){\circle*{1}}
  \put(35,17.5){\circle*{1}}
  \put(30,40){\circle*{1}}
  \put(35,47.5){\circle*{1}}
  \put(45,47.5){\circle*{1}}
    \end{picture} &
     \setlength{\unitlength}{1mm}
    \begin{picture}(50,60)(0,5)
   \put(10,10){\line(2,-1){5}}
   \put(15,7){\line(1,0){10}}
   \put(15,7){\line(1,3){10}}
  \put(10,10){\line(1,1){10}}
  \put(10,10){\line(1,3){10}}
  \put(25,7){\line(4,1){12}}
  \put(20,20){\line(1,0){15}}
  \put(30,50){\line(1,0){15}}
  \put(20,40){\line(2,-1){5}}
  \put(20,40){\line(1,1){10}}
  \put(25,37){\line(1,0){10}}
  \put(35,37){\line(4,1){12}}
  \put(25,7){\line(1,3){10}}
  \put(20,20){\line(1,3){10}}
  \put(37,10){\line(1,3){10}}
  \put(37,10){\line(1,2){2.5}}
  \put(47,40){\line(1,2){2.5}}
  \put(40,15){\line(1,3){10}}
  \put(45,50){\line(1,-1){5}}
  \put(35,20){\line(1,3){10}}
  \put(35,20){\line(1,-1){5}}
  \put(3,2){\makebox(10,10){$0$}}
  \put(17,10){\makebox(10,10){$A$}}
  \put(5,20){\makebox(15,10){$I$}}
  \put(10,10){\circle*{1}}
  \put(15,7){\circle*{1}}
  \put(25,7){\circle*{1}}
  \put(20,40){\circle*{1}}
  \put(20,20){\circle*{1}}
  \put(37,10){\circle*{1}}
  \put(35,20){\circle*{1}}
  \put(25,37){\circle*{1}}
  \put(30,50){\circle*{1}}
  \put(45,50){\circle*{1}}
  \put(40,15){\circle*{1}}
  \put(50,45){\circle*{1}}
  \put(47,40){\circle*{1}}
  \put(35,37){\circle*{1}}
      \end{picture}\\
  \mbox{(a) $P=J_1+J_2+I$} & \mbox{(b) $P=A+I$}
 \end{tabular}
 \end{center}
\caption{A polygonal prism $P$}
 \label{fig1}
\end{figure}
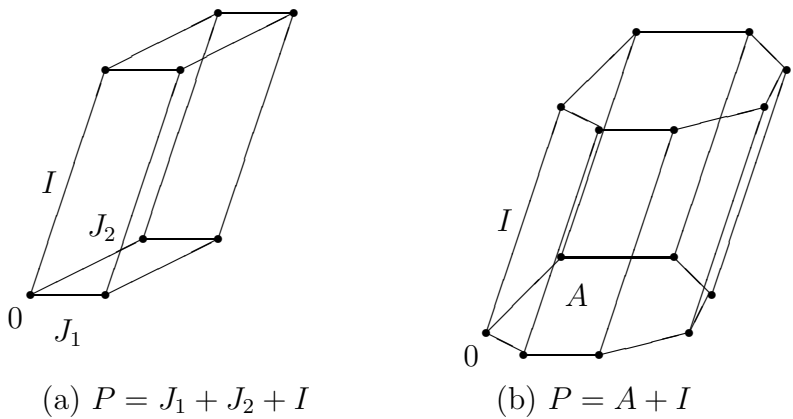

In this section we investigate the normality of a lattice polytope $P$ of dimension three
which is the Minkowski sum of a lattice polygon $A$ of dimension two and a lattice
line segment $I$.
See Figure~\ref{fig1}~(b).
Here we set 
$M=M'\oplus L$, $\mbox{rank}\  M'=2$, $\mbox{rank}\  L=1$ and 
$A\subset M'_{\mathbb{R}}$.

If $A$ is a parallelogram, that is, if $A$ is the Minkowski sum $J_1+J_2$ of
two not parallel lattice line segments  $J_1$ and $J_2$, 
then $P=A+I$ is normal because it is a parallelotope.
See Figure~\ref{fig1}~(a).

From this observation we obtain the following proposition.

\begin{prop}\label{2;p1}
Let $A$ be a nonsingular lattice polygon in $M'_{\mathbb{R}}$ not isomorphic to 
a basic triangle and $I$ a lattice line segment not contained in 
$M'_{\mathbb{R}}$.
Then $P=A+I$ is normal.

In the case that  $A$ is a basic triangle, if $I\subset L_{\mathbb{R}}$, then 
$P=A+I$is normal.
Here a lattice triangle is called basic if it is isomorphic to 
the convex hull of the origin and a basis
of $M'\cong \mathbb{Z}^2$.
\end{prop}

In order to prove Proposition~\ref{2;p1}, it is enough to show the following lemma.

\begin{lem}
If a nonsingular lattice polygon $A\subset M'_{\mathbb{R}}$ is not basic,
then it is covered by a union of lattice parallelograms.
\end{lem}
{\it Proof.}
Take a coordinates $(x,y)$ in $M'_{\mathbb{R}}$.
Assume that a lattice polygon $A$ is $r$-gonal, that is, $A$ has $r$ edges.
By a suitable affine transformation of $M'$, we may take a vertex $v_0$ of $A$ to be
the origin, an edge $E_1$ from $v_0$ to be on the positive part of the $x$-axis and the other
edge $E_r$ from $v_0$ on the positive part of the $y$-axis. 
We will prove the lemma by dividing into several steps.

(a): The case of $A=\mbox{Conv}\{0, (a,0), (0,1), (b,1)\}$, that is, $r=4$.  
If $a=b$, then $A$ is a parallelogram.
Set $a<b$ and $s=b-a$.  If we set $A_i=\mbox{Conv}\{0, (a,0), (i,1), (a+i,1)\}$ for 
$i=0,1, \dots, s$, then $A_i$ is a parallelogram and $A$ is covered by the union of $A_i$ 
with $i=0, \dots, s$.  The same is when $a>b$.

(b):  Set $F(E_1):=A\cap(0\le y\le1)$.  Since $A$ is nonsingular, 
$F(E_1)$ is also a lattice polygon.  From (a), we see that $F(E_1)$ is covered by
a union of lattice parallelograms.
For all edges $E_1, \dots, E_r$ of $A$, define $F(E_i)$ in the same way.
Then we have
$$
A = A^\circ \cup \bigcup_{i=1}^r F(E_i),
$$
where $A^\circ$ is the convex hull of $\mbox{Int}(A)\cap M'$.

If $\dim A^\circ\le1$, then $A$ is covered by the union of $F(E_i)$.
If $\dim A^\circ=2$, then $A^\circ$ is a nonsingular lattice polygon.
If $A^\circ$ is not isomorphic to a basic triangle, then we continue this process.

(c):  When $A^\circ$ is isomorphic to a basic triangle, we may consider $A$ is the $4$ times
multiple $\mbox{Conv}\{0, (4,0),(0,4)\}$ of a basic triangle, or,
a polygon obtained from this by cutting off several basic triangles at vertices.
Set $A':=A\cap(1\le y\le4)$.  Then we have a decomposition $A=A'\cup F(E_1)$
and we see that $A'$ is nonsingular and $\dim (A')^\circ\le1$.
Thus we see that $A$ is covered by a union of lattice parallelograms
in this case.

Since the normalized area of $A^\circ$ is an integer less than that of $A$, 
this process is stop after
several steps.
\hfill  $\Box$

\bigskip

Next we introduce another direct sum decomposition of $M=M'\oplus L$ with
respect to $I$ of the Minkowski sum $P=A+I$.

Set $L':=(\mathbb{R}I)\cap M\cong\mathbb{Z}$ and $M= M''\oplus L'$ with
the projection map $\pi:M\to M''$.
We note that $M'=(\mathbb{R}A)\cap M\cong \mathbb{Z}^2$ does not always 
coincide with $M''$.
Set $B:=\pi(A)\subset M''_{\mathbb{R}}$.  Then $B$ is a lattice polygon in $M''_{\mathbb{R}}$.

\begin{figure}[h]
 \begin{center}
    \setlength{\unitlength}{1mm}
    \begin{picture}(50,60)(0,10)
   \put(10,10){\line(2,-1){5}}
   \put(15,7){\line(1,0){10}}
   \put(15,7){\line(0,1){10}}
  \put(10,10){\line(1,1){10}}
  \put(10,10){\line(0,1){10}}
  \put(25,7){\line(4,1){12}}
  \put(20,20){\line(1,0){15}}
  \put(10,20){\line(1,2){10}}
  \put(10,20){\line(2,-1){5}}
  \put(15,17){\line(2,1){10}}
  \put(25,22){\line(1,1){12}}
  \put(37,34){\line(1,3){3}}
  \put(25,7){\line(0,1){15}}
  \put(20,20){\line(0,1){20}}
  \put(37,10){\line(0,1){24}}
  \put(37,10){\line(1,2){2.5}}
  \put(40,43){\line(-2,5){5}}
  \put(40,15){\line(0,1){28}}
  \put(20,40){\line(1,1){15}}
  \put(35,20){\line(0,1){35}}
  \put(35,20){\line(1,-1){5}}
  \put(3,2){\makebox(10,10){$v_0$}}
  \put(25,8){\makebox(10,10){$B$}}
   \put(23,35){\makebox(10,10){$A$}}
  \put(0,10){\makebox(15,10){$I$}}
  \put(10,10){\circle*{1}}
  \put(15,7){\circle*{1}}
  \put(25,7){\circle*{1}}
  \put(10,20){\circle*{1}}
  \put(20,20){\circle*{1}}
  \put(37,10){\circle*{1}}
  \put(35,20){\circle*{1}}
  \put(15,17){\circle*{1}}
  \put(20,40){\circle*{1}}
  \put(25,22){\circle*{1}}
  \put(40,15){\circle*{1}}
  \put(35,55){\circle*{1}}
  \put(37,34){\circle*{1}}
  \put(40,43){\circle*{1}}
      \end{picture}
   \end{center}
\caption{An upright polygonal prism $Q(A)$}
 \label{fig2}
\end{figure}
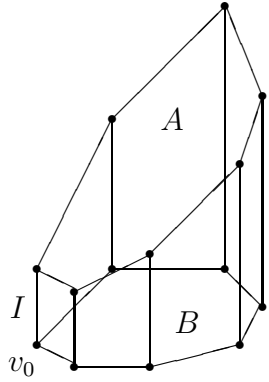

Take  coordinates $(x,y)$ in $M''_{\mathbb{R}}$ and $z$ in $L'_{\mathbb{R}}$.
From a suitable affine transform of $M$, we may set so that a vertex $v_0$ of $P=A+I$
is the origin and $P$ is contained in the upper half space 
$(z\ge0)$.  Then define $Q(A)$ as the convex hull of $B\times 0$ and $A+I$.
The polytope $Q(A)$ is an upright polygonal prism with  the $r$-gonal polygon $B$ as
its base and $A$ as its roof.  See Figure~\ref{fig2}.

\begin{prop}\label{2;p2}
Let $A$ be a nonsingular lattice polygon in $M'_{\mathbb{R}}$
not isomorphic to a basic triangle.
Then $Q(A)$ defined above is normal.
\end{prop}
{\it Proof.} 
Decompose the lattice polygon $B$ into a union of basic lattice triangles
$B_i$
($i=1, \dots, s$) with vertices in $B\cap M''$.
For each $B_i$, define $R(B_i)$ as the convex hull of $(B_i\times\mathbb{R}_{\ge0})
\cap Q(A)\cap M$.
Then the prism $R(B_i)$ is normal.  We have a cover 
$$
Q(A) = (A+I)\cup \bigcup_{i=1}^s R(B_i).
$$
Since $A+I$ is normal from Proposition~\ref{2;p1}, the  polytope $Q(A)$ is normal.
\hfill  $\Box$

\bigskip

\section{Union of polygonal prisms}

In this section we assume that a projective toric fibered 3-fold $X$ over $\mathbb{P}^1$
has one irreducible invariant fiber.

As in Section 1, we set $N_0:=\mathbb{Z}$ and $\Delta_0:=\{\mathbb{R}_{\le0}, \{0\}, \mathbb{R}_{\ge0}\}$.   Then $X(\Delta_0)=\mathbb{P}^1$.
The torus equivariant morphism $\varphi: X=X(\Delta) 
\to \mathbb{P}^1$ is defined by the morphism
of fans $\varphi^{\sharp}: (N, \Delta) \to (N_0, \Delta_0)$ with 
$\varphi^{\sharp}(N)=N_0$.
Set $N_0^{\vee}$ the dual to $N_0$.   Denote by $L$ the image of the dual homomorphism
$\varphi*: N_0^{\vee} \to M$.  Then we have a direct sum decomposition
$M=M_f\oplus L$, where $M_f^{\vee}\cong N_f:=(\varphi^{\sharp})^{-1}(0)$.
The subset $\Delta_f:=\{\sigma\in \Delta; \varphi^{\sharp}(\sigma)=0\}$ is a fan of $N_f$.
A general fiber of $\varphi$ is the toric surface $X(\Delta_f)$.

Let $\mathcal{L}$ be an ample line bundle on a toric fibered 3-fold $X(\Delta)$ over
$\mathbb{P}^1$.  Let $P$ be the lattice polytope in $M_{\mathbb{R}}$ coresponding to
the polarized toric 3-fold $(X(\Delta), \mathcal{L})$.
Denote by $\mathcal{L}_f$ the restriction of the ample line bundle $\mathcal{L}$ to 
$X(\Delta_f)$.   The polarized toric surface 
$(X(\Delta_f), \mathcal{L}_f)$ defines a nonsingular lattice polygon 
$B\subset (M_f)_{\mathbb{R}}$.  Then the lattice polytope 
$P$ is contained in the polygonal prism $B\times \mathbb{R} \subset (M_f \oplus L)_{\mathbb{R}}
=M_{\mathbb{R}}$ and each side wall of the prism 
contains a facet of $P$.

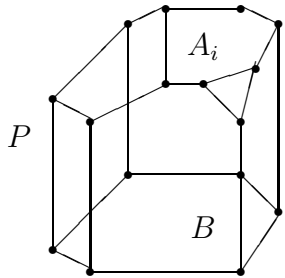
\begin{figure}[h]
 \begin{center}
    \setlength{\unitlength}{1mm}
    \begin{picture}(50,40)(0,10)
   \put(10,10){\line(2,-1){5}}
   \put(15,7){\line(1,0){20}}
   \put(15,7){\line(0,1){20}}
  \put(10,10){\line(1,1){10}}
  \put(10,10){\line(0,1){20}}
  \put(35,7){\line(2,3){5}}
  \put(20,20){\line(1,0){15}}
  \put(10,30){\line(1,1){10}}
  \put(10,30){\line(2,-1){5}}
  \put(15,27){\line(2,1){10}}
  \put(25,32){\line(1,0){5}}
  \put(25,32){\line(0,1){10}}
  \put(37,34){\line(1,2){3}}
  \put(30,32){\line(1,-1){5}}
  \put(30,32){\line(3,1){7}}
  \put(20,20){\line(0,1){20}}
  \put(35,7){\line(0,1){20}}
  \put(25,42){\line(1,0){10}}
  \put(40,40){\line(-2,1){5}}
  \put(40,15){\line(0,1){25}}
  \put(20,40){\line(2,1){5}}
  \put(35,27){\line(1,4){2}}
  \put(35,20){\line(1,-1){5}}
  %\put(3,2){\makebox(10,10){$v_0$}}
  \put(25,8){\makebox(10,10){$B$}}
   \put(25,32){\makebox(10,10){$A_i$}}
  \put(-2,20){\makebox(15,10){$P$}}
  \put(10,10){\circle*{1}}
  \put(15,7){\circle*{1}}
  \put(35,7){\circle*{1}}
  \put(10,30){\circle*{1}}
  \put(20,20){\circle*{1}}
  \put(25,42){\circle*{1}}
  \put(35,20){\circle*{1}}
  \put(15,27){\circle*{1}}
  \put(20,40){\circle*{1}}
  \put(25,32){\circle*{1}}
    \put(30,32){\circle*{1}}
  \put(40,15){\circle*{1}}
  \put(35,42){\circle*{1}}
   \put(35,27){\circle*{1}}
  \put(37,34){\circle*{1}}
  \put(40,40){\circle*{1}}
      \end{picture}
   \end{center}
\caption{Union of upright polygonal prisms}
 \label{fig3}
\end{figure}

If one invariant fiber of $\varphi$ is irreducible, then $P$ has a facet isomorphic to $B$.
We may draw the picture of $P$ so that it is a polygonal upright prism with $B$ as the
base and the roof consists of a collection of lattice polygons.
See Figure~\ref{fig3}.

\medskip

\begin{prop}\label{3;p1}  Assume that a nonsingular projective 
toric fibered 3-fold $\varphi: X(\Delta)\to \mathbb{P}^1$ has one irreducible
invariant fiber $\varphi^{-1}([1:0])$, that is, one irreducible invariant fiber is isomorphic
to a general fiber.   Then an ample line bundle on 
$X(\Delta)$ is always normally generated.
\end{prop}
{\it Proof.}
Take coordinates $(x,y,z)$ in 
$M_{\mathbb{R}}=(M_f\oplus L)_{\mathbb{R}}$ so that $(M_f)_{\mathbb{R}}=(z=0)$.
Let $P\subset M_{\mathbb{R}}$ be a lattice polytope corresponding to an ample
line bundle $\mathcal{L}$ on $X(\Delta)$.

From our assumption, $P$ has the special facet $B$ corresponding to the irreducible
 fiber of $\varphi^{-1}([1:0])$.
From a suitable affine transform of $M$, we may assume that $P$ is contained 
in the upper half space $(z\ge0)$ and $B$ is contained in the plane $(z=0)$.

Set $A_1, \dots, A_s$ the all facets in the roof of the upright prism $P$.
Set $B_i=\pi(A_i)$ the lattice polygon in $(M_f)_{\mathbb{R}}$
defined as the image of a facet $A_i$ by the projection 
$\pi: (M_f\oplus L)_{\mathbb{R}}\to (M_f)_{\mathbb{R}}$.
For each facet $A_i$ define 
$Q(A_i):=(B_i\times L_{\mathbb{R}})\cap P$.
Then we have a decomposition of $P$ as a union of polygonal prisms $Q(A_i)$.
For each 
$A_i$, set $M_i:=(\mathbb{R}A_i)\cap M \cong \mathbb{Z}^2$.
Then 
$A_i$ is a nonsingular lattice polygon in $(M_i)_{\mathbb{R}}$.

If $A_i$ is not a basic triangle, then $Q(A_i)$ is normal from Proposition~\ref{2;p2}.
Even if $A_i$ is a basic triangle if it meets a side wall of $P$, then it is normal because 
$M_i\oplus L\cong M$.

We assume that $A_i$ is a basic triangle and meets no side walls of $P$.
Set $v_1,v_2,v_3$ the three vertices of $A_i$ and $E_1,E_2,E_3$
the edges of $P$ from $v_1, v_2, v_3$ outside $A_i$, respectively.
Let $w_j$ be the lattice points on the edge $E_j$ nearest $v_j$ for $j=1,2,3$.
Set $\tilde{A_i}:=\mbox{Conv}\{w_1, w_2, w_3\}$.
Then the lattice triangle $\tilde{A_i}$ is similar and parallel to $A_i$
since $P$ is nonsingular.
If $\tilde{A_i}\cong A_i$, then $P=Q(A_i)$.  It contradicts the assumption.
Thus $\tilde{A_i}$ is not basic. 
The subset $(\pi(\tilde{A_i})\times L_{\mathbb{R}})\cap P$ of $P$ can be decomposed into 
a union of the slice $\mbox{Conv}\{A_i, \tilde{A_i}\}$ of the roof and the rest 
$Q(\tilde{A_i})$.  Both are normal.

Since $P$ is covered by a union of normal lattice polytopes, it is normal.  \hfill  $\Box$

\bigskip

\section{General case}

\begin{prop}\label{4;p1}
Let $X(\Delta)$ be a projective nonsingular toric fibered 3-fold over $\mathbb{P}^1$.
Then an ample line bundle on $X(\Delta)$ is always normally generated.
\end{prop}

In this section we assume that two invariant fibers of $\varphi: X(\Delta) \to \mathbb{P}^1$
are reducible.

Let $\cal{L}$ be an ample line bundle on $X(\Delta)$ and $P\subset M_{\mathbb{R}}$
the lattice polytope corresponding to $(X(\Delta),\cal{L})$.
As in the proof of Proposition~\ref{3;p1}, take coordinates $(x,y,z)$ in
$M_{\mathbb{R}}=(M_f\oplus L)_{\mathbb{R}}$ so that $(M_f)_{\mathbb{R}}=(z=0)$.
From a suitable affine transform of $M$, we may assume that $P$ is contained 
in the upper half space $(z\ge0)$.

Let $F$ be a side wall, a facet of $P$ parallel to $L$.  
From a suitable affine transform of $M$, if we set $F$ to be contained in the plane $(y=0)$,
then we define the lattice polytope $P(F):=P\cap(0\le y\le1)$.
Then $P(F)$ is normal by Lemma 2.5 in \cite{Od}.
Set $P^\circ:=\mbox{Conv}\{\mbox{Int}(P)\cap M\}$ the inner polytope of $P$.
If $\dim P^\circ\le 2$, then $P$ is normal because $P$ is covered by a union of
$P(F)$ for all side walls $F$ of $P$.

(S1) \ We assume $\dim P^\circ=3$.

By the projection $\pi: (M_f\oplus L)_{\mathbb{R}} \to (M_f)_{\mathbb{R}}$,
we define $G:=\pi(P)$, which is a nonsingular lattice polygon.
Set $G^\circ:=\mbox{Conv}\{\mbox{Int}(G)\cap M_f\}$ the inner polygon of $G$,
which is also a nonsingular lattice polygon from the assumption (S1).

Set $A_1, \dots, A_s$ the all facets in the roof of the upright prism $P$ 
 for $i=1,\dots,s$.
For a lattice point $m\in M$, denote by $l(m)$ the line through $m$ parallel to 
the $z$-axis.

\begin{lem}\label{lem:2}
For a lattice point $m'\in (\partial G^\circ)\cap M_f$, the line segment $l((m',0))\cap P$
has length greater than or equal to two.
\end{lem}
{\it Proof}. We note that both edges of $l((m',0))\cap P$ are lattice points.
If the length is less than two, then it is one, hence $P$ is contained in $(0\le z\le1)$.
In this case, two invariant fibers of $\varphi$ are both irreducible. This contradicts to the assumption 
in the beginning of this section.  \hfill $\Box$

\bigskip

Set $A_i':=A_i\cap\pi^{-1}(G^\circ)$ for $i=1, \dots, s$.
Set $I:=[0,(0,0,-1)]$ the unit interval on the $z$-axis of negative direction.
Define the convex sets as 
$$
U_1(P):=\mbox{Conv}\{\cup_i (A_i'+I)\} \quad \mbox{and} \quad
U_2(P):=\mbox{Conv}\{\cup_i (A_i'+2I)\}.
$$
From Lemma~\ref{lem:2}, we have $U_1(P)\subset U_2(P)\subset P$.

In order to prove that $P\cap M+P\cap M =(2P)\cap M$,
take a lattice point $m$ in $2P$ and consider $\frac12m\in P$.
If $m$ is located on the boundary of $2P$, then it is a lattice point on a 
lattice polytope of dimension less than three, hence, there exist two lattice points
$m_1,m_2\in (\partial P)\cap M$ such that $m=m_1+m_2$.

We may assume that $\frac12m$ is contained in the interior of $P$ and that $\frac12m\notin M$.
If $\frac12m$ is contained in $P(F)$ for a side wall $F$ of $P$, then there exist $m_1, m_2\in P(F)$
such that $m=m_1+m_2$ because $P(F)$ is normal.
Here we call that $m$ is {\it 2-normal in} $P$ if there exist $m_1,m_2\in  P\cap M$ 
such that $m=m_1+m_2$.

\begin{lem}\label{lem:3}
If $\frac12m$ is contained in $U_1(P)$, then $m$ is 2-normal.
\end{lem}
{\it Proof}.
If $A_i'$ is not a basic triangle, then $I+A_i'$ is normal.  Evev if $A_i'$ is basic, 
as in the proof of Proposition~\ref{3;p1} take $\widetilde{A_i'}$ a triangle similar and parallel to
$A_i'$ so that $I+\widetilde{A_i'}$ is normal.
Thus if $\frac12m$ is contained in a $I+A_i'$, then $m$ is 2-normal in $P$.

Let $\{m_1', \dots, m_r'\}$ be the set of vertices of $G^\circ$
Set $\widetilde{m_i}\in \pi^{-1}(G^\circ)\cap M$ the lattice point on the roof of $U_1(P)$.
For three $\widetilde{m_i},\widetilde{m_j},\widetilde{m_k}$ of $\{\widetilde{m_1}, \dots, \widetilde{m_r}\}$,
set $F_{ijk}=\mbox{Conv}\{\widetilde{m_i},\widetilde{m_j},\widetilde{m_k}\}$.
Decompose the lattice triangle $\pi(F_{ijk})$ into a union of basic lattice triangle $G_l$, $l=1, \dots, t$
with vertices in $M_f$.
For $G_l$, define $R(G_l)$ the convex hull of $(G_l\times L)\cap P\cap M$.  Then the prism 
$R(G_l)$ is normal. $F_{ijk}$ is covered bu a union of $R(G_l)$'s.
We take all three of $\{\widetilde{m_1}, \dots, \widetilde{m_r}\}$.
Then the boundary of $U_1(P)$ is covered by a union of normal polytopes.

Assume that $\frac12m$ is not contained in $I+A_i'$ nor $R(G_l)$.
If $\pi(\frac12m)\in M_f$, then $\frac12m$ is contained in a lattice line segment parallel to $L$, hence,
$m$ is 2-normal.

Assume that  $\pi(\frac12m)\notin M_f$.
Since $G^\circ=\cup_i \pi(A_i')$, we can choose a $\pi(A_i')$ such that $\pi(\frac12m)\in \pi(A_i')$.
Decompose $\pi(A_i')$ into a union of basic lattice triangle $G_j$ with vertices in $M_f$
such that $G_1=\mbox{Conv}\{u_1,u_2,u_3\}$ with $\frac12m =\frac12(u_1+u_2)$
since $2\pi(\frac12m)\in M_f$.

Set
$$
\mbox{Conv}\{l((u_1,0))\cap P \cap M\}=[u_1^-,u_1^+] \quad \mbox{and} \quad 
\mbox{Conv}\{l((u_1,0))\cap P \cap M\}=[u_2^-, u_2^+].
$$
Here the $z$-coordinates of $u_i^+$ is bigger than those of $u_i^-$.
Since $u_1^+$ and $u_2^+$ are contained in $U_1(P)$, the lattice line segment $[u_1^+,u_2^+]$
is contained in $U_1(P)$.

Consider the lattice line segment $[u_1^-,u_2^-]$.
Im $(M_f)_{\mathbb{R}}$, the line containing both $u_1$ and $u_2$ meets opposite edges of
$G^\circ$. From two edges we take vertices $\{m_1',m_2',m_3',m_4'\}$ such that
the line segment $[u_1,u_2]$ is contained in the interior of the lattice polygon
$\mbox{Conv}\{m_1',m_2',m_3',m_4'\}$.
Set $m_1'',m_2'',m_3'',m_4''$ the vertices of $U_1(P)$ on the bottom with $\pi(m_i'')=m_i'$.
Consider the convex hull of $(\mathbb{R}_{\ge0}I +\mbox{Conv}\{m_1'',m_2'',m_3'',m_4''\})\cap P$.
This contains $u_1^-$ and $u_2^-$, hence, the line segment $[u_1^-,u_2^-]$ because $U_2(P)\subset P$.

Consider the lattice polygon $R=\mbox{Conv}\{u_1^-,u_1^+, u_2^-,u_2^+\}\subset P$.
Since $R$  contains $\frac12m$, we see that $m$ is 2-normal in $P$.
\hfill $\Box$

\bigskip

{\it Proof of Proposition~\ref{4;p1}}.
Set $B_1, \dots, B_t$ be all facets in the bottom of $P$. Set $B_j':=B_j\cap \pi^{-1}(G^o)$.
We define 
$$
D_1(P):=\mbox{Conv}\{\cup_j(B_j'+(-I))\}\subset P.
$$
As in the proof of Lemma~\ref{lem:3}, we see that if $\frac12m$ is contained in $D_1(P)$,
then $m$ is 2-normal.

Assume that $\frac12m$ is not contained in $U_1(P)$ nor $D_1(P)$.
As in the proof of Lemma~\ref{lem:3}, we can choose a basic lattice triangle
$\mbox{Conv}\{u_1,u_2,u_3\}\subset \pi(A_i')$ with $frac12m=\frac12(u_1+u_2)$.
In this case, since $u_1^-, u_2^-\in D_1(P)$, the lattice line segment $[u_1^-,u_2^-]$ is contained in
$D_1(P)$. Thus $m$ is 2-normal.
\hfill  $\Box$

\bigskip

\noindent
{\bf Remark}.
Nonsingularity condition of a 
toric fibered 3-fold $\varphi: X(\Delta) \to \mathbb{P}^1$ in Proposition~\ref{4;p1} is necessary.
We know a singular toric fibered 3-fold over $\mathbb{P}^1$ with 
 a very ample but not
normally generated line bundle on it.
Finally, we will give an example found by Burns and Gubeladze \cite[Exercise 2.24]{BG}.

For a positive integer $q$, define a lattice tetrahedron as 
$$
Q_q:=\mbox{Conv}\{0, (1,0,0), (0,1,0), (1,1,q)\}.
$$
If $q\ge2$, then $Q_q$ is not very ample.
Set $I=[0, (0,0,1)]$ the unit interval on the $z$-axsis.
Define $P_q:=Q_q+I$ sa the Minkowski sum.
Let $(X, \mathcal{L})$ be the polarized toric 3-fold corresponding to $P_q$.
Then this $X$ is a singular toric fibered 3-fold over $\mathbb{P}^1$
and $\mathcal{L}$ is very ample.
If $q\ge 4$, then $\mathcal{L}$ is not normally generated.  See also \cite{Og2}.

\end{document}